\documentclass{article}
\newtheorem{theorem}{Theorem}
\newtheorem{rem}{Remark}

\newtheorem{cor}{Corollary}
\newtheorem{exmp}{Example}
\newtheorem{prop}{Proposition}

\renewcommand{\epsilon}{\varepsilon}

\renewcommand{\hat}{\widehat}

\usepackage{epsfig}
\title{Bounds for contractive semigroups and second order systems}
\author{Kre\v simir Veseli\'c\thanks{Fernuniversit\" at Hagen, Lehrgebiet
Mathematische Physik, Postfach 940, D-58084 Hagen, Germany, e-mail:
kresimir.veselic@fernuni-hagen.de.}}
\date{ }
\begin{document}
\maketitle
\begin{abstract} We derive a uniform bound for the difference of two
contractive semigroups, if the difference of their 
generators is form-bounded by the  
Hermitian parts of the generators themselves.
We construct a semigroup dynamics for second order systems
with fairly general operator coefficients and apply our bound to
the perturbation of the damping term. The result is illustrated
on a dissipative wave equation. As a consequence the exponential
decay of some second order systems is proved. 
\end{abstract}

\section{Introduction}
\label{sec:introduction}
The aim of this paper is to derive a new perturbation bound
for strongly continuous contractive
semigroups in a Hilbert space and to apply it to damped systems
of second order.
Let \(e^{At}\), \(e^{Bt}\) be strongly continuous contractive
semigroups in a Hilbert space
\({\cal X}\). Their generators are maximal dissipative in the sense of
\cite{Kato} (as negatives of accretive operators, introduced there; we
will follow the notations and the terminology of this monograph).
 
We consider a rather restricted kind of perturbation, it reads
formally
\begin{equation}\label{relboundformal}
|(x,(B - A)y)|^2 \leq \varepsilon^2
\Re(-Bx,x)\Re(-Ay,y),\quad \varepsilon > 0.
\end{equation}
As a result we obtain a uniform estimate for the semigroups:
\[
\|e^{Bt} - e^{At}\| \leq
\frac{\varepsilon}{2}.
\]
Note that here we have not the classical situation: 'unperturbed
object plus a small perturbation' in which the perturbed object
often has first to be constructed and then the distance between the two
is measured (see e.g. \cite{Kato} Ch. XI, Th.~2.1). 
We impose no condition whatsoever
on the size of the positive constant \(\varepsilon\) but we know that
both \(A\) and \(B\) are dissipative, and both operators
appear in a symmetric way.
Moreover,
no requirements are made about the size of the subspace
\({\cal D}(A)\cap{\cal D}(B)\), it could even be trivial. To this end,
(\ref{relboundformal}) is rewritten
in a 'weak form' as 
\[
|(B^*x,y) - (x,Ay)|^2 \leq \varepsilon^2
\Re(-B^*x,x)Re(-Ay,y).
\] 

This kind of perturbation
will appear to be the proper setting for treating semigroups,
generated by second order systems
\begin{equation}\label{MCK}
M\ddot{x} + C\dot{x} + Kx = 0.
\end{equation}
Here \(M\), \(C\), \(K\) can be finite symmetric matrices,
with {\em the mass matrix} \(M\) positive semidefinite, 
{\em the stiffness matrix}
\(K\) positive definite and {\em the damping matrix} \(C\) 
accretive\footnote{For simplicity we use the term 'damping matrix'
for \(C\) although it is not necessarily symmetric and thus 
may include a gyroscopic component.}
(our result seems to be new even in the matrix case).
Or, \(M\), \(C\), \(K\) may  be differential operators with
similar properties. We will construct a 
contractive semigroup, naturally attached to (\ref{MCK}),
where \(M\), \(C\), \(K\) are
understood as sesquilinear forms satisfying some mild natural
regularity conditions. This 
construction seems to cover damped systems, more general than those
treated in previous literature (cf. e.g. \cite{HryShk1}, \cite{HryShk2},
\cite{NSV}), for instance, \(M\) is allowed to have 
a nontrivial null-space and \(C\) need not be symmetric. 
Subsequently,
we derive a semigroup bound for such systems 
in which the damping term \(C\) is subject to a perturbation
of the same type as (\ref{relboundformal}). 
As a consequence, the exponential decay of some
damped systems will be proved. In particular, under the 
additional assumption that 
\(C\) be sectorial, a second order system
is exponentially stable, 
if and only if the system with the 'pure symmetric damping'
\(\hat{C} = (C^* + C)/2\) is such. 
In these applications an important 
property of the condition (\ref{relboundformal}) will be used:
it is invariant under the inversion of both operators.\footnote{
A related perturbation result for finite
matrices was proved in \cite{Vesb}}

The article is organised as follows. In Sect.~2 we prove the main result
in a 'local' and a 'global' version.
We also include an
analogous bound for discrete semigroups, although we have no
application for it as yet.

In Sect.~2 we apply this theory to abstract damped systems of the form
(\ref{MCK}), including the construction of the semigroup itself.
In Sect.~3 we apply our theory to the damped wave equation in one
dimension. 

\section{An abstract perturbation bound}
Let \(A\) be the generator of a strongly continuous semigroup in a 
Hilbert space \({\cal X}\). By \({\cal T}(A)\) we denote the set of 
all differentiable semigroup trajectories
\[
S = \left\{x = e^{At}x_0, t \geq 0\right\},
\mbox{ for some }x_0 \in {\cal D}(A).
\]
\begin{theorem}\label{pert1} Let \(A\), \(B\) be the generators of 
strongly continuous semigroups in a
Hilbert space \({\cal X}\) (then \(A^*\), \(B^*\) are also such).
Suppose that there exist trajectories \(S \in {\cal T}(A)\),
\(T \in {\cal T}(B^*)\) and an \(\varepsilon > 0\) such that 
for any \(y \in S\), \(x \in T\) 
\begin{equation}\label{relbound}
|(B^*x,y) - (x,Ay)|^2 \leq \varepsilon^2
\Re(-B^*x,x)\Re(-Ay,y).
\end{equation} 
Then for all such \(x,y\)
\begin{equation}\label{weakestimate}
|\left(x,(e^{Bt} - e^{At})y)\right)|
\leq \frac{\varepsilon}{2}\|x\|\|y\|.
\end{equation}
(Note that in (\ref{relbound}) it is tacitly assumed that 
the factors on the right hand side are non-negative.)
\end{theorem}
{\bf Proof.} For \(y \in S\), \(x \in T\) we have
\[
\frac{d}{ds}\left(e^{B^*s}x,e^{A(t - s)}y\right) =
\]
\[
\left(e^{B^*s}B^*x,e^{A(t - s)}y\right) -
\left(e^{B^*s}x,e^{A(t - s)}Ay\right),
\]
which is continuous in \(s\), so by integrating from \(0\)
to \(t\) we obtain the weak Duhamel formula
\[
(e^{B^*t}x,y) - (x,e^{At}y) =
\int_0^t\left[(B^*e^{B^*t}x,e^{A(t-s)}y) -
(e^{B^*t}x,Ae^{A(t-s)}y)\right]ds.
\]
By using (\ref{relbound}) and the Cauchy-Schwarz inequality it 
follows
\[
|\left(x,(e^{Bt} - e^{At}y)\right)|^2 \leq
\]
\[
\left(\int_0^t|(B^*e^{B^*t}x,e^{A(t-s)}y) -
(e^{B^*t}x,Ae^{A(t-s)}y)|ds\right)^2 \leq
\]
\[
\varepsilon^2
\left(\int_0^t\sqrt{\Re(-B^*e^{B^*s}x,e^{B^*s}x)
\Re(-Ae^{A(t-s)}y,e^{A(t-s)}y)}ds\right)^2 \leq
\]
\[
\varepsilon^2
\int_0^t\Re(-B^*e^{B^*s}x,e^{B^*s}x)ds
\int_0^t\Re(-Ae^{As}y,e^{As}y)ds.
\]
By partial integration we compute 
\[
{\cal I}(A,y,t) = \int_0^t\Re(-Ae^{As}y,e^{As}y)ds =
\]
\[
-\|e^{As}y\|^2\Big|_0^t - {\cal I}(A,y,t),
\]
\begin{equation}\label{IAt}
{\cal I}(A,y,t) = \frac{1}{2}\left((y,y) - \|e^{At}y\|^2\right).
\end{equation}
Obviously
\[
0 \leq {\cal I}(A,y,t) \leq \frac{1}{2}(y,y)
\]
and \({\cal I}(A,y,t)\) increases with \(t\). Thus, there exist 
limits
\[
{\cal I}(A,y,t) \nearrow {\cal I}(A,y,\infty) =
\frac{1}{2}\left(y,y) - P(A,y)\right), \quad t \to \infty
\]
\[
\|e^{At}y\|^2 \searrow P(A,y), \quad t \to \infty
\]
with
\[
0 \leq {\cal I}(A,y,\infty) \leq \frac{1}{2}(y,y). 
\]
(and similarly for \(B^*\)). Altogether
\begin{equation}\label{fineestimate}
|\left(x,(e^{Bt} - e^{At}y)\right)|^2\leq
\frac{\varepsilon^2}{4}
\left((x,x) - P(B^*,x)\right)\left((y,y) - P(A,y)\right)
\end{equation}
\[
\leq \frac{\varepsilon^2}{4}(x,x)(y,y).
\]
Q.E.D.\\

\begin{rem}
{\em As a matter of fact, in the proof above neither of the operators
need be densely defined. In this case the assertion of the theorem 
is valid only in the weak form 
\[
|(e^{B^*t}x,y) - (x,e^{At}y)| 
\leq \frac{\varepsilon}{2}\|x\|\|y\|.
\]
}
\end{rem}
\begin{cor}\label{expdecay}
Suppose that (\ref{relbound}) holds for all \(y\) from some
\(S \in {\cal T}(A)\) and all \(x \in {\cal D}\), where \({\cal D}\)
is a dense subspace, invariant under \(e^{B^*t},\ t \geq 0\). Then
\begin{equation}\label{normestimate}
\|\left(e^{Bt} - e^{At}y\right)\| \leq
\frac{\varepsilon}{2}\|y\|.
\end{equation}
\end{cor}
By setting \(\varepsilon = 0\) in (\ref{normestimate}) we obtain the 
known uniqueness of the solution of a first order differential
equation:
\[
e^{Bt}y = e^{At}y,\quad t \geq 0.
\]

If both \(A\) and \(B\) are maximal dissipative then
\begin{equation}\label{PA}
P(A,y) = (P(A)y,y), \quad
P(A) = \mbox{s-}\lim_{t\to\infty}e^{A^*t}e^{At}.
\end{equation}
The strong limit \(P(A)\) above exists by the contractivity and obviously
\( 0 \leq P(A) \leq I\) in the sense of forms (and similarly for \(B^*\)).
\begin{cor} \label{unifestimatecor}
If (\ref{relbound}) holds for all \(x \in {\cal D}\),
\(y \in {\cal E}\), where \({\cal D}\), \({\cal E}\) are dense subspaces,
invariant under \(e^{B^*t}\), \(e^{At}\), respectively, then
\begin{equation}\label{weakfine}
|\left(x,(e^{Bt} - e^{At})y)\right)|^2
\leq \frac{\varepsilon^2}{4}\left((x,x) - (P(B^*)x,x))\right)
\left((y,y) - (P(A)y,y)\right).
\end{equation}
In particular,
\begin{equation}\label{uniformestimate}
\|e^{Bt} - e^{At}\| \leq
\frac{\varepsilon}{2}.
\end{equation}
\end{cor}
\begin{rem}{\em The corollary above certainly holds, if (\ref{relbound})
is fulfilled for all \(x \in {\cal D}(B^*)\) and all 
\(y \in {\cal D(A)}\)
(it is enough to require the validity of (\ref{relbound}) on respective cores)
and this will be the situation in our applications. In any of these cases
both \(B^*\) and \(A\) (and then also \(B\) and \(A^*\)) are, in fact,
maximal dissipative. 
}
\end{rem}
The condition (\ref{relbound}) has a remarkable property of being
{\em inversion invariant} i.e.~\(B^*\) and \(A\) may be replaced
by their inverses.
\begin{prop} \label{inverses} Suppose that 
both \(B^*\) and \(A\) (and then also 
\(B\) and \(A^*\)) are (not necessarily boundedly) invertible. 
Then (\ref{relbound}),
valid for all \(x \in {\cal D}(B^*)\) and all
\(y \in {\cal D(A)}\) is equivalent to
\begin{equation}\label{relboundinv}
|(B^{-*}\xi,\eta) - (\xi,A^{-1}\eta)|^2 \leq
\varepsilon^2\Re(-B^{-*}\xi,\xi)\Re(-A^{-1}\eta,\eta),\quad
\xi \in {\cal D}(B^{-*}), \quad \eta \in {\cal D}(A^{-1}).
\end{equation}
\end{prop}
{\bf Proof.} Just set \(B^*x = \xi,\ Ay = \eta\). Q.E.D.\\

Note that in all our results above no further restriction to the
constant \(\varepsilon\) was imposed. This is partly due to the 
fact that the perturbation is measured by both the ``perturbed''
and the ``unperturbed'' operator in a completely symmetric way.
This kind of perturbation bound will prove particularly appropriate
for our applications below. If \(\varepsilon\) is further restricted
important new conclusions can be drawn.

A semigroup is called {\em exponentially stable}\footnote{Some
authors call this property the {\it uniform} exponential
stability.} or {\em exponentially decaying}, if
\begin{equation}\label{exp_stab}
\|e^{At}\| \leq ce^{-\beta t},\quad t\geq 0
\end{equation}
for some \(c, \beta > 0\).

\begin{cor} \label{decayAB}
If in Corollary \ref{unifestimatecor} we have \(\varepsilon < 2\) 
then the exponential decay of one of the semigroups implies the same
for the other.
\end{cor}
{\bf Proof.} Just recall that the exponential stability follows,
if \(\|e^{At}\| < 1 \) for some \(t > 0\). Q.E.D.
\begin{rem}
{\em
In all that was said thus far there is an obvious symmetry: 
in (\ref{relbound}) we may replace
\(A,\ B^*\) by \(B,\ A^*\), thus obtaining the dual estimate
\begin{equation}\label{relbound_dual}
|(Bx,y) - (x,A^*y)|^2 \leq \varepsilon^2
\Re(-Bx,x)\Re(-A^*y,y).
\end{equation}
with completely analogous results. Obviously, (\ref{relbound}) and
(\ref{relbound_dual}) are equivalent, if \({\cal D}({\cal A} =
{\cal D}({\cal A}^*)\) and \({\cal D}({\cal B} =
{\cal D}({\cal B}^*)\). 
}
\end{rem}

{\bf Discrete semigroups.} Every step of the perturbation theory,
developed above can be correspondingly extended to discrete semigroups.
An operator \(T\) is called a contraction, if
\(\|T\| \leq 1\). For any such operator \(T\) the strong limit 
\[
Q(T) = \mbox{s-}\lim_{n\to\infty}T^{*n}T^n
\]
obviously exists and satisfies
\[
0 \leq Q(T) \leq I.
\]
The following theorem sums up the most important facts.
\begin{theorem}\label{discrete_theorem} Let \(A\), \(B\)
be contractions and
\begin{equation}\label{relbound_discr}
|((B - A)x,y)|^2 \leq \varepsilon^2 
((I - B^*B)x,x)((I - AA^*)y,y)
\end{equation}
for all \(x,y\)  and some \(\varepsilon \geq 0\) (note that in
(\ref{relbound_discr}) the right hand side is always non-negative). Then
\begin{equation}\label{weakestimate_diskr_xy}
|((B^n - A^n)x,y)|^2 \leq
\varepsilon^2
((I - Q(B))x,x)((I - Q(A^*))y,y)
\end{equation}
and, in particular,
\begin{equation}\label{weakestimate_diskr}
\|B^n - A^n\|
\leq \varepsilon
\sqrt{\|I - Q(A^*)\|\|I - Q(B)\|}
\leq \varepsilon.
\end{equation}
\end{theorem}
{\bf Proof.} For any \(x,y\) we have
\[
|((B^n - A^n)x,y)|^2 = |(\sum_{k=0}^{n-1}A^k(B - A)B^{n-k-1}x,y)|^2 \leq
\]
\[
\left(\sum_{k=0}^{n-1}|((B - A)B^{n-k-1}x,A^{*k}y)|\right)^2
\]
\[
\leq \varepsilon^2 
\left(\sum_{k=0}^{n-1}\sqrt{(A^k(I - AA^*)A^{*k}y,y)
(B^{*n-k-1}(I - B^*B)B^{n-k-1}x,x)}\right)^2
\]
\[
\leq \varepsilon^2
\sum_{k=0}^{n-1}(A^k(I - AA^*)A^{*k}y,y)
\sum_{k=0}^{n-1}(B^{*k}(I - B^*B)B^ky,y)
\]
\[
= \varepsilon^2
((I - A^nA^{*n})y,y)((I - B^{*n}B^n)x,x)
\]
and (\ref{weakestimate_diskr_xy}) follows.
Here we have used the identity
\begin{equation}\label{dentity_diskr_AB}
\sum_{n=0}^{n-1}A^k(I - AB)B^k = I - A^nB^n.
\end{equation}
Q.E.D.\\

It may be interesting to note that (\ref{dentity_diskr_AB}) 
appears to be a discrete analog of
\begin{equation}\label{dentity_cont_AB}
\int_0^te^{A\tau}(A + B)e^{B\tau}d\tau =
-\left(I -e^{At}e^{Bt} \right)
\end{equation}
on which (\ref{IAt}) was based.

Any contraction \(A\) is exponentially stable, if and only if
\(\|A^n\| < 1\) for some \(n\). This leads to a result, analogous to
Cor. \ref{decayAB}.
\begin{cor}\label{decayAB_discr}
Let \(A\) and \(B\) be contractions satisfying (\ref{weakestimate_diskr_xy})
with \(\varepsilon < 1\). Then the exponential stability of
one of them implies the same for the other.
\end{cor}
One might wonder that the bound (\ref{uniformestimate})
is uniform in \(t\) although the involved semigroups need not be 
exponentially decaying. As a simple example consider dissipative
operators \(A \), \(B\) in a finite dimensional space. Then 
each of these operators is known to be an orthogonal sum
of a skew-Hermitian part and an exponentially stable part. By
(\ref{relbound}) (which is now equivalent to (\ref{relboundformal}))
the skew-Hermitian parts of \(A \) and \(B\) coincide and the
difference \(e^{Bt} - e^{At}\) decays exponentially. The situation with 
discrete semigroups is similar.

In the infinite dimensional case 
the uniformity of the bound (\ref{uniformestimate})
is a more serious fact as will be illustrated on applications
from Mathematical Physics below.

\section{Application to damped systems}
An abstract damped linear system is governed by a formal second order
differential equation in a vector space \({\cal Y}_0\)
\begin{equation}\label{MCKforms}
\mu(\ddot{y},v) + \theta(\dot{y},v) + \kappa(y,v) = 0,
\end{equation}
where \(\mu\), \(\theta\), \(\kappa\) are sesquilinear forms
with the following properties:
\begin{itemize}
\item \(\kappa\) symmetric, strictly positive,
\item \(\mu\) symmetric, positive, \(\kappa\)-closable,
\item \(\theta\) \(\kappa\)-bounded, accretive.
\end{itemize}
A possible way to turn (\ref{MCKforms}) into an operator equation
is to take \((u,v)=\kappa(u,v)\) as the scalar product and to
complete accordingly \({\cal Y}_0\) to a Hilbert space \({\cal Y}\).
By the known representation theorems (\cite{Kato}) we have
\begin{equation}\label{generateMC0}
\mu(u,v) = (Mu,v),\quad \theta(u,v) = (Cu,v),
\end{equation}
where \(M\) is (possibly unbounded) selfadjoint and positive
and \(C\) is bounded accretive. We now replace (\ref{MCKforms}) by
\begin{equation}\label{MCequation}
M\ddot{y} + C\dot{y} + y = 0,
\end{equation}
where the time derivatives \(\dot{y}\), \(\ddot{y}\) are taken in 
\({\cal Y}\).\footnote{Our choice of the underlying scalar product
in \({\cal Y}_0\) is fairly natural but not the only relevant one.
One could show that very different, even topologically non-equivalent,
choices of the scalar product still lead to the essentially 
same semigroup dynamics, see \cite{NSV}.}

To the equation (\ref{MCequation}) one naturally associates the
phase space system, obtained by the formal substitution
\begin{equation}\label{substitution}
x_1 = y,\quad x_2 = M^{1/2}\dot{y}
\end{equation}
which leads to the first order equation
\[
\frac{d}{dt}
\left(\begin{array}{c}
x_1 \\
x_2 \\
\end{array}\right)=
{\cal A}
\left(\begin{array}{c}
x_1 \\
x_2 \\
\end{array}\right)
\]
with
\begin{equation}\label{Aformal}
{\cal A} =
\left(\begin{array}{rr}
0          &    M^{-1/2}        \\
-M^{-1/2}  &  -M^{-1/2}CM^{-1/2}\\
\end{array}\right) .
\end{equation}
which then should generate a contractive semigroup which realises the
dynamics. Our conditions are far too general for this
\({\cal A}\) to make sense as it stays (note that \(M\) may
have a nontrivial null-space). However, the formal inverse
\begin{equation}\label{Aplus}
{\cal A}^+ =
\left(\begin{array}{rr}
-C         &   -M^{1/2}        \\
M^{1/2}    &        0          \\
\end{array}\right)
\end{equation}
is more regular, although not necessarily bounded. Considered in the
'total energy' Hilbert space
\(\hat{{\cal X}} = {\cal Y} \oplus {\cal Y}\), \({\cal A}^+\) has
the following properties
\begin{equation}\label{property1}
{\cal A}^+ \mbox{ is maximal dissipative,}
\end{equation}
\begin{equation}\label{property2}
{\cal D}({\cal A}^+) = {\cal D}(({\cal A}^+)^*)
= {\cal D}(M^{1/2}) \oplus {\cal D}(M^{1/2}), 
\end{equation}
\begin{equation}\label{property3}
{\cal N}({\cal A}^+) = {\cal N}(({\cal A}^+)^*).
\end{equation}
All this follows from the fact that \({\cal A}^+\) is
a sum of the skew-selfadjoint operator
\[
\left(\begin{array}{rr}
 0         &   -M^{1/2}        \\
M^{1/2}    &        0          \\
\end{array}\right)
\]
and a bounded dissipative operator
\[
-
\left(\begin{array}{rr}
 C         &        0          \\
 0         &        0          \\
\end{array}\right).
\]
Thus, \({\cal N}({\cal A}^+)\) reduces both \({\cal A}^+\) and its adjoint,
the same is the case with
\begin{equation}\label{X}
{\cal X} = {\cal N}({\cal A}^+)^\bot.
\end{equation}
More precisely, \({\cal A}^+\) is a direct sum of the null operator
and a maximal dissipative invertible operator \({\cal A}^{-1}\) 
in the Hilbert space
\({\cal X}\), defined on
\[
{\cal X}\cap{\cal D}({\cal A}^+)
\]
which is dense in \({\cal X}\). Obviously, the operator \({\cal A}\) 
is again maximal dissipative
and this is by definition the generator of our
semigroup. The space \({\cal X}\) may be called {\em the physical
phase space} for the system (\ref{MCequation}).\footnote{
A different but related construction was used
in \cite{NSV} where
both \(M\) and \(C\) are symmetric, but possibly unbounded.}

Denoting by \(Q\) the orthogonal projection onto the space
\({\cal X}\) in \(\hat{{\cal X}}\) we have, in fact,
\begin{equation}\label{pseudor}
(\lambda - {\cal A})^{-1}Q = 
\frac{1}{\lambda} - \frac{1}{\lambda^2}
\left(\frac{1}{\lambda} - {\cal A}^+\right)^{-1},
\quad \Re\lambda \neq 0,
\end{equation}
which is immediately verified.
\begin{prop} The null-space \({\cal N}({\cal A^+})\) 
satisfies the inclusion
\begin{equation}\label{NAplusincl}
{\cal N}({\cal A^+})\supseteq 
\left({\cal N}(C) \cap {\cal N}(M)\right) \oplus {\cal N}(M).
\end{equation}
If, in addition, \(C\) is sectorial then we have the equality
\begin{equation}\label{NAplus}
{\cal N}({\cal A^+}) =
\left({\cal N}(C) \cap {\cal N}(M)\right) \oplus {\cal N}(M).
\end{equation}
\end{prop}
{\bf Proof.}\footnote{In the case of \(C\) symmetric and \(M\)
bounded this formula was proved in \cite{NSV}.} 
Now, \({\cal N}({\cal A^+})\) is given by the equations
\[
-Cx_1 - M^{1/2}x_2 = 0,\quad M^{1/2}x_1 = 0, \quad 
x_{1,2} \in {\cal D}(M^{1/2}).
\]
From this the inclusion (\ref{NAplusincl}) follows.
Let now \(C\) be sectorial. The above equations
imply \((Cx_1,x_1) = -(M^{1/2}x_2,x_1) = 0\). By the 
assumed sectoriality it follows \(Cx_1 = 0\), so
(\ref{NAplus}) follows. Q.E.D.\\

The fact that the semigroup dynamics exists only on a closed subspace
\({\cal X}\) of \(\hat{{\cal X}}\) is quite natural, even in the 
finite dimensional space: one cannot prescribe velocity initial data
on the parts of the space where the mass is vanishing. If \(M\)
is injective --- no matter how singular \(M^{-1}\) may be ---
our dynamics exists on the whole space \(\hat{{\cal X}}\).

It can be shown
(\cite{NSV}) that this semigroup provides an appropriate solution to the 
second order system (\ref{MCequation}) via the formulae
(\ref{substitution}), at least in the special case of
\(M,\ C\) bounded symmetric. In our, more general situation
we can show that \({\cal A}\) yields the ``true'' dynamics
by way of approximation. We approximate the operator \(M\) by a sequence 
\(M_n\) of bounded, positive operators such that
\begin{equation}\label{Mapprox}
M_n^{1/2}x \to M^{1/2}x, \quad x \in {\cal D}(M^{1/2}).
\end{equation}
If, in addition, all \(M_n\) are positive definite the operator
(\ref{Aformal}) 
\begin{equation}\label{invconv}
{\cal A}_n =
\left(\begin{array}{rr}
0               &    M_n^{-1/2}             \\
-M_n^{-1/2}  &  -M_n^{-1/2}CM_n^{-1/2}\\
\end{array}\right) 
\end{equation}
is bounded dissipative in \(\hat{{\cal X}}\) and its semigroup
trivially reproduces the solution of the so modified second order system
\begin{equation}\label{MCequation_n}
M_n\ddot{y} + C\dot{y} + y = 0,
\end{equation}
An example of such sequence is
\[
M_n = f_n(M),\quad
f_n(\lambda) =
\left\{\begin{array}{rl}
\frac{1}{n}, & 0 \leq \lambda \leq \frac{1}{n}\\
\lambda,     &  \frac{1}{n} \leq \lambda \leq n\\
n,           & n \leq \lambda
\end{array}\right.
\]
Note that here, in addition, the operators \(M_n\) are both
bounded and boundedly invertible, being positive definite.
\begin{prop}For any 
\(x =
\left(\begin{array}{c}
x_1 \\
x_2 \\
\end{array}\right)
\in {\cal X}\) 
and any approximation sequence
(\ref{Mapprox}) we have
\begin{equation}\label{convergence}
e^{{\cal A}_n}x \to e^{{\cal A}}x, \quad n \to \infty. 
\end{equation}
uniformly on any compact interval in \(t\). Choose, in addition, 
\(M_n\) as positive definite and set
\[
\left(\begin{array}{c}
y_n(t) \\
u_n(t) \\
\end{array}\right)
= e^{{\cal A}_nt}
\left(\begin{array}{c}
x_1 \\
x_2 \\
\end{array}\right),\quad
\left(\begin{array}{c}
y(t)\\
u(t)\\
\end{array}\right)
= e^{{\cal A}t}
\left(\begin{array}{c}
x_1 \\
x_2 \\
\end{array}\right).
\]
Then \(y_n(t)\) solves (\ref{MCequation_n}) with 
\(u_n(t) = M_n^{1/2}\dot{y}_n(t)\) and
\[
y_n(t) \to y(t),\quad M_n^{1/2}\dot{y}_n(t) \to M^{1/2}\dot{y}(t),
\quad n \to \infty.
\] 
\end{prop}
{\bf Proof.} By (\ref{Mapprox}) we have \({\cal A}_n^{-1} \to {\cal A}^+\)
in the strong resolvent sense (see \cite{Kato}, Ch.~VIII, Th,~1.5)
i.e.
\[
(\lambda - {\cal A}_n^{-1})^{-1} \to (\lambda - {\cal A}^+)^{-1},
\quad \Im \lambda \neq 0.
\]
Hence by (\ref{pseudor}),
\[
(\lambda - {\cal A}_n)^{-1} =
\frac{1}{\lambda} - \frac{1}{\lambda^2}
\left(\frac{1}{\lambda} - {\cal A}_n^{-1}\right)^{-1}
\to
\]
\[
\frac{1}{\lambda} - \frac{1}{\lambda^2}
\left(\frac{1}{\lambda} - {\cal A}^+\right)^{-1}
= (\lambda - {\cal A})^{-1}Q.
\]
all in the strong sense. Now the Trotter-Kato convergence theory 
(\cite{Kato}) can be applied to give 
\begin{equation}\label{Trotter-Kato}
e^{{\cal A}_\eta t}x \to e^{{\cal A}t}x ,\quad \eta \to 0
\end{equation}
for all \(x \in {\cal X}\). (The original Trotter-Kato
theorem requires the injectivity of the strong limit
in (\ref{invconv}), but the same proof
is easily seen to accomodate our slightly more general setting.)
The remaining assertions are now straightforward.
Q.E.D.\\

We now apply our abstract theory from Sect.~2 to a 
second order system with variable damping. 
\begin{theorem}\label{pertA}
Let
\[
{\cal A}^+ =
\left(\begin{array}{rr}
-C         &   -M^{1/2}        \\
M^{1/2}    &        0          \\
\end{array}\right), \quad
\hat{{\cal A}}^+ =
\left(\begin{array}{rr}
-\hat{C}         &   -M^{1/2}        \\
M^{1/2}          &        0          \\
\end{array}\right)
\]
where \(M\) is bounded, positive selfadjoint and \(C\), \(\hat{C}\) are 
bounded accretive operators satisfying
\begin{equation}\label{Cbound}
\big|\left((\hat{C} - C)x,y\right)\big|^2 \leq 
\varepsilon^2\Re(Cy,y)\Re(\hat{C}x,x)
\end{equation}
for all \(x,y \in {\cal Y}\) and some \(\varepsilon > 0\).
Then \({\cal A}^+\) and  \(\hat{{\cal A}}^+\) have the same null-space
and the respective contractive semigroup generators \({\cal A}\) and 
\(\hat{{\cal A}}\) 
in \({\cal X}\) from (\ref{X}) satisfy the assumptions of 
Cor.~\ref{unifestimatecor}, in particular, 
\begin{equation}\label{uniformestimatecal}
\|e^{\hat{{\cal A}}t} - e^{{\cal A}t}\| \leq
\frac{\varepsilon}{2}.
\end{equation}
\end{theorem}
{\bf Proof.} Obviously (\ref{Cbound}) is equivalent
to
\begin{equation}\label{Abound}
\Big|\left((\hat{{\cal A}}^+ - {\cal A}^+)x,y\right)\Big|^2 \leq
\varepsilon^2\Re(-{\cal A}^+y,y)\Re(-\hat{{\cal A}}^+x,x)
\end{equation}
for all \(x,y \in {\cal D}({\cal A}^+) = {\cal D}(\hat{{\cal A}}^+)\). 
From this it follows that
 \({\cal A}^+\) and \(\hat{{\cal A}}^+\) 
have the same null-space. Furthermore, by (\ref{property2}) 
the domains of the four operators \({\cal A}^+, \hat{{\cal A}}^+,
({\cal A}^+)^*, (\hat{{\cal A}}^+)^*\)
coincide and (\ref{Abound}) is equivalent to both
(\ref{relbound}) and (\ref{relbound_dual}) for 
for \(A = {\cal A}^+\) and \(B = \hat{{\cal A}}^+\) and then also
for \(A = {\cal A}^{-1}\) and \(B = \hat{{\cal A}}^{-1}\)
(note that in our situation we have \(\Re(-{\cal A}^+y,y) =
\Re(-({\cal A}^+)^*y,y)\) and \(\Re(-\hat{{\cal A}}^+x,x) =
\Re(-(\hat{{\cal A}}^+)^*x,x)\)). 
Now apply Prop.~\ref{inverses}
 and Cor.~\ref{unifestimatecor}. Q.E.D.\\

Note the important role of the 'inverse-invariance property'
in Prop.~\ref{inverses} in the proof above because we have 
no explicit formulae for the generators \({\cal A}\) and
\(\hat{{\cal A}}\) and there is no control on their domains of definition. 
\begin{rem}{\em
It is important to note that our perturbation bound (\ref{Cbound})
can be readily expressed in the language of the original forms
in (\ref{MCKforms}). Obviously, (\ref{Cbound}) is equivalent to
\begin{equation}\label{formbound}
\big|(\hat{\gamma} - \gamma)(x,y)\big|^2 \leq
\varepsilon^2\Re\gamma(y,y)\Re\hat{\gamma}(x,x).
\end{equation}
}
\end{rem}

We now prove some stability results for second
order systems.
\begin{theorem}\label{pertA_ineq} Let the system (\ref{MCequation})
be exponentially stable\footnote{By the exponential stablity
of a second order system we mean the exponential stablity of
the generated semigroup.} with a symmetric \(C = C^{(1)}\) and let
\[
0 \leq C^{(1)} \leq D \leq \alpha C^{(1)}.
\]
Then the exponential stability holds with \(C = D\) and vice versa.
\end{theorem}
{\bf Proof.} Set
\[
C_k = C^{(1)} + \frac{k}{n}(D - C^{(1)}),\quad k = 0,\ldots,n.
\]
Then \(C_0 = C^{(1)}\), \(C_n = D\) and
\[
0 \leq C_{k+1} - C_k \leq \frac{\alpha - 1}{n}C^{(1)}
\]
and
\[
|((C_{k+1} - C_k)x,y)|^2 \leq
((C_{k+1} - C_k)x,x)((C_{k+1} - C_k)y,y)
\]
\[
\leq \left(\frac{\alpha - 1}{n}\right)^2(C^{(1)}x,x)(C^{(1)}y,y)
\leq \left(\frac{\alpha - 1}{n}\right)^2(C_{k+1}x,x)(C_ky,y).
\]
Now choose \(n > (\alpha - 1)/2\) and use Theorem \ref{pertA} 
and Corollary
\ref{decayAB}. Use induction: the exponential stability carries 
over from \(C_k\) to \(C_{k+1}\) and vice versa. Q.E.D.\\

In particular, the exponential stability with \(C\)
implies the same with \(\alpha C\) for any positive \(\alpha\).
A similar technique can be applied to gyroscopic systems:
\begin{theorem}\label{pertA_sector}
Suppose that in (\ref{MCequation}) the operator \(C\) 
is sectorial. Then the exponential stability of this system
is equivalent to the exponential stability of the 'purely damped'
system
\begin{equation}\label{MChatequation}
M\ddot{y} + \hat{C}\dot{y} + y = 0,\mbox{ with }
\hat{C} = \frac{C^* + C}{2}.
\end{equation}
\end{theorem}
{\bf Proof.} By sectoriality there exists \(N > 0\) such that
\begin{equation}\label{ReImC}
|\Im(Cy,y)| \leq N\Re(Cy,y) \mbox{ for all } x,y. 
\end{equation}
We have
\[
\Re(Cx,x) = (\hat{C}x,x), \quad
-i\Im(Cx,x) = -i((\hat{C} - C)x,x).
\]
The operators \(\hat{C}\) and \(-i(\hat{C} - C)\) are symmetric, so 
the inequality (\ref{ReImC}) may be polarised to read
\[
((\hat{C} - C)x,y)|^2 \leq N^2\Re(Cy,y)\Re(Cx,x) 
= N^2\Re(Cy,y)\Re(\hat{C}x,x).
\]
Assume first \(N < 2\). Then apply Theorem \ref{pertA} and Corollary
\ref{decayAB} to obtain the exponential stability with
\(C\). Now drop the condition \(N < 2\) and proceed by
induction. Introduce the sequence
\[
C_k = \hat{C} + \frac{k}{n}(C - \hat{C}),\quad k = 0,\ldots,n.
\]
Then obviously
\[
|((C_{k+1} - C_k)x,y)|^2 = \frac{1}{n^2}|((\hat{C} - C)x,y)|^2\leq
\]
\[
\frac{N^2}{n^2}\Re(Cy,y)\Re(\hat{C}x,x) =
\frac{N^2}{n^2}\Re(C_ky,y)\Re(C_{k+1}x,x).
\]
Now choose \(n < 2/N\) and apply the above consideration 
to the consecutive pairs \(C_k\), \(C_{k+1}\). We may begin at the bottom
with \(C_0 = \hat{C}\) or at the top with \(C_n = C\).
Q.E.D.\\

The perturbations allowed in our theory are smooth enough not
to change the null-space of \({\cal A}^+\) and thus the physical
phase space remains the same. More general perturbations will
possibly violate this property, and the resulting approximations will,
like (\ref{Trotter-Kato}), hold only on subspaces.
%\[
%\eta = \frac{\varepsilon}{2\sqrt{1 - \varepsilon}}
%\]
%we see that to the value \(\eta = 1\) there corresponds
%\(\varepsilon = \varepsilon_0 = \sqrt{8} - 2 \approx 0.83\). Thus,
%we have
\section{The damped wave equation}
Here we apply our general theory to the wave equation in one dimension
\begin{equation}\label{wave}
\rho(x)u_{tt} + \gamma(x)u_t - (d(x)u_{tx})_x - (k(x)u_x)_x = 0
\end{equation}
for the unknown function \(u = u(x,t)\), \(a < x < b\) and 
\(0< t < \infty\).
The functions \(\rho(x),\ \gamma(x),\ d(x),\ k(x)\) are assumed
to be non-negative and measurable; in addition, \(\rho(x),\ \gamma(x)\)
are bounded and
\begin{equation}\label{akbound}
\mbox{ess inf}_{a < x < b}k(x) > 0,\quad \mbox{ess sup}_{a < x < b}
\frac{d(x)}{k(x)} < \infty.
\end{equation}
 The boudary conditions are
\begin{equation}\label{boundary1}
u(a,t) = 0, \quad u_x(b,t) + \zeta u_t(b,t) = 0,\quad \zeta \geq 0.
\end{equation}
%Here the cases \(a = -\infty\) or \(b = \infty\) are allowed;
%then the corresponding boundary condition is dropped.

This is a formally dissipative equation which we shall understand in
its weak form
\begin{equation}\label{waveweak}
\mu(u_{tt},v) + \theta(u_t,v) + \kappa(u,v) = 0
\end{equation}
with \(u(a) = v(a) = 0\) and
\begin{equation}\label{mu}
\mu(u,v) = \int_a^b\rho(x)u\bar{v}dx,
\end{equation}
\begin{equation}\label{theta}
\theta(u,v) = \int_a^b\left(\gamma(x)u\bar{v} + 
d(x)u'\bar{v}'\right)dx + \zeta u(a)\bar{v}(b),
\end{equation}
\begin{equation}\label{kappa}
\kappa(u,v) = \int_a^b k(x)u'\bar{v}'dx.
\end{equation}
The forms \(\mu,\ \theta\) are symmetric and positive.  \(\theta\) is obviously 
\(\kappa\)-bounded while  \(\mu\) is \(\kappa\)-closable. 
As the underlying Hilbert space 
\({\cal Y}\) we take the functions
with the scalar product 
\begin{equation}\label{Yproduct}
(u,v) =
\kappa(u,v) = \int_a^b k(x)u'\bar{v}'dx,\quad u(a) = v(a) = 0.
\end{equation}
Then under our conditions,
\begin{equation}\label{generateMC}
\mu(u,v) = (Mu,v),\quad \theta(u,v) = (Cu,v)
\end{equation}
where \(M,\ C\) are positive selfadjoint operators, with
bounded \(C\) and \(M\). 
Thus, we end up
with the second order system (\ref{MCequation}) and
(\ref{waveweak}) gives rise to a contractive semigroup on  the space
\({\cal X}\) which is determined from the null-spaces
of \(M,\ C\). 

Note that in order for \(M\) to have a non-trivial
null-space it is not sufficient that the function \(\rho\) vanishes
just on a set of positive measure, rather \(\rho\) must vanish 
on an interval (and similarly for \(C\)). If \(\rho\) vanishes on an interval
and \(\gamma\) does not, then (\ref{wave}) is of mixed type
(hyperbolic - parabolic). All such cases are covered by our theory. 

Now for the perturbation. We perturb the damping parameters
\(\gamma(x),\ d(x),\ \zeta\) into 
\(\hat{\gamma}(x),\ \hat{d}(x),\ \hat{\zeta}(x)\), which satisfy the same
conditions as \(\gamma(x),\ d(x),\ \zeta\) above and are such that
\begin{equation}\label{pertgamma}
|\hat{\gamma}(x) - \gamma(x)| \leq \varepsilon
\sqrt{\hat{\gamma}(x) \gamma(x)}
\end{equation}
\begin{equation}\label{perta}
|\hat{d}(x) - d(x)| \leq \varepsilon
\sqrt{\hat{d}(x) d(x)}
\end{equation}
\begin{equation}\label{pertzeta}
|\hat{\zeta} - \zeta| \leq \varepsilon
\sqrt{\hat{\zeta} \zeta}
\end{equation}
This is a 'relatively small' change of the damping parameters,
commonly encountered in practice. The corresponding operators \(C\)
and \(\hat{C}\) are immediately seen to satisfy
(\ref{Cbound}). Hence Theorem \ref{pertA} applies and the corresponding
semigroups satisfy (\ref{uniformestimatecal}).

One might be interested to obtain perturbation results under the more
common assumptions involving only the 'unperturbed' data and the perturbation:
\begin{equation}\label{pertgamma0}
|\hat{\gamma}(x) - \gamma(x)| \leq \eta 
\gamma(x)
\end{equation}
\begin{equation}\label{perta0}
|\hat{d}(x) - d(x)| \leq \eta
d(x)
\end{equation}
\begin{equation}\label{pertzeta0}
|\hat{\zeta} - \zeta| \leq \eta
\zeta
\end{equation}
with \(\eta < 1\). This implies (\ref{pertgamma}) -- (\ref{pertzeta})
with 
\[
 \varepsilon = \frac{\eta}{\sqrt{1 - \eta}}.
\]
But the real use of (\ref{pertgamma0}) -- (\ref{pertzeta0})
consists merely in insuring the non-negativity of the perturbed
damping parameters and the conditions (\ref{akbound});
all this is usually known in advance, so there is
no need to abandon the much less restrictive conditions
(\ref{pertgamma}) -- (\ref{pertzeta}).

In view of Corollary \ref{expdecay} we conclude that
{\em if the equation (\ref{wave}) decays exponentially with the 
damping parameters \(\gamma(x),\ d(x),\ \zeta\), then the same 
will be the case with \(\hat{\gamma}(x),\ \hat{d}(x),\ \hat{\zeta}\),
if the constant \(\varepsilon\) is less than \(2\)}.  Theorem
\ref{pertA_ineq} also applies accordingly.
The situation in higher dimensions is similar and the results are
analogous.\\

As a second example consider the equation (\ref{wave}) on the 
infinite interval \(0 < x < \infty\) with the boundary condition
\begin{equation}\label{boundary2}
u(0,t) = 0.
\end{equation}
For simplicity we take
\begin{equation}\label{krho1}
k(x) = \rho(x) \equiv 1,
\end{equation}
whereas \(\gamma(x) \geq 0\) is supposed to satisfy
\begin{equation}\label{gammabound}
D = \sup_{u \in {\cal Y}}
\frac{\int_0^\infty\gamma(x)|u(x)|^2dx}{\int_0^\infty|u'(x)|^2dx}
< \infty,
\end{equation}
where \({\cal Y}\) is the set of all \(u\) which are absolutely
continuous, vanish at zero and have a square integrable \(u'\); 
this is obviously a Hilbert space with the scalar product
\[
(u,v) = \int_0^\infty u(x)\bar{v}'(x)dx.
\]
The class of functions \(\gamma\) satisfying (\ref{gammabound})
is not void since it includes
\[
\gamma(x) = \frac{1}{x^2},\mbox{ with } D = 4
\]
(\cite{Kato} Ch.~VI 4.1).
The form
\[
\mu(u,v) = \int_0^\infty u(x)\bar{v}(x)dx 
\]
defined on \({\cal D}(\mu) = L^2(0,\infty)\cap{\cal Y}\)
is closed in \({\cal Y}\), so (\ref{generateMC}) yields
a positive unbounded operator \(M\) with a trivial null-space and
a bounded \(C\). Hence our semigroup construction applies
and under the perturbation (\ref{pertgamma}) the bounds 
(\ref{uniformestimate}), (\ref{weakfine}) hold. 

Such semigroups are in general not exponentially 
decaying (they are usually extended to uniformly
bounded groups) and
will give rise to a non-trivial scattering theory on an
'absorbing obstacle' represented by the short range damping
function \(\gamma\). Further considerations along these lines go beyond
the scope of this article.
%A further detailed investigation of the solutions of this
%dissipative wave equation, although interesting and relevant,
%goes beyond the scope of this note.

\end{document}